\documentclass{elsart}
\usepackage{amssymb, amsfonts}
\usepackage{amsmath}

\def\A{{\mathcal A}} 
\def\R{{\mathbb R}}
\let\bt\beta
\let\dlt\delta
\let\eps\varepsilon
\let\w\omega
\let\nowt\varnothing
\def\sbs{\subseteq}
\def\ol{\overline}
\newcommand{\mdn}{\textrm{m}\dlt\textrm{n}}
\newcommand{\sq}{\subseteq}
\newcommand{\sm}{\smallsetminus}
\newcommand{\les}{\leqslant}
\newcommand{\n}{\varnothing}

\begin{document}
\begin{frontmatter}
\title{Monotone Versions of Countable Paracompactness}

\author{Chris Good},
\thanks{c.good@bham.ac.uk}
\author{Lylah Haynes}
\thanks{haynesl@maths.bham.ac.uk}
\address{School of Mathematics, University of Birmingham, Birmingham B15
2TT, UK}

\begin{abstract} One possible natural monotone version of countable
paracompactness,
MCP, turns out to have some interesting properties.
We investigate various other possible monotonizations of
countable paracompactness and how they are related.
\end{abstract}

\begin{keyword}
Monotone countable paracompactness \sep MCP \sep monotone $\delta$-normality\\
\emph{AMS subject classification:} 54E20, 54E30
\end{keyword}
\end{frontmatter}

\section{Introduction}

Countably paracompact spaces were introduced by Dowker \cite{dowk}
and Kat\v etov \cite{kat} and their importance lies in Dowker's
classic result: $X\times[0,1]$ is normal if and only if $X$ is both
normal and countably paracompact. Countable paracompactness occurs
with normality in a number of other results concerning both
separation and the insertion of continuous functions (see
\cite{eng,gks}). Moreover, it turns out that a number of
set-theoretic results concerning the separation of closed discrete
collections in normal spaces have direct analogues for countably
paracompact spaces. For example, Burke \cite{bu} modifies Nyikos's
\lq provisional\rq\ solution to the normal Moore space problem by
showing that countably paracompact, Moore spaces are metrizable
assuming PMEA and Watson \cite{wat} shows that, assuming V=L, first
countable, countably paracompact spaces are collectionwise
Hausdorff. When normality is strengthened to monotone normality,
pathology is reduced and the need for set-theory in such results is
generally avoided. It turns out that the same is often true for one
possible monotone version of countable paracompactness, MCP (see
below for a definition), introduced in \cite{gks} and also in
\cite{txl} and \cite{pan}. For example, every MCP Moore space is metrizable and
every first countable or locally compact MCP space is collectionwise
Hausdorff. Interestingly, although every monotonically normal space
is collectionwise Hausdorff, if there is an MCP space which fails to
be collectionwise Hausdorff, then there is a measurable cardinal,
and if there are two measurable cardinals, then there is an MCP
space that is not collectionwise Hausdorff \cite{gk}.

Since there is a reasonable amount one can say about MCP spaces and since
there is a large number of characterizations of countable paracompactness, it
makes
sense to ask about other monotone versions of countable paracompactness. In this
paper we look at some of these possible properties.
Our notation is standard as used in \cite{eng} and \cite{hand},
where any undefined terms may be found. All spaces are $T_1$ and
regular.

The characterizations of countable paracompactness fall, broadly,
into four categories relating to covering properties, perfect
normality (see conditions (2-4) of Theorem \ref{cp equivs}),
interpolation theorems (conditions (5-7)), and normality in products
(8-10). It seems harder to make sensible, useful and widely
satisfied definitions of monotone covering properties and so we do
not consider countable paracompactness as a covering property in
this paper; neither do we consider $\aleph_0$-expandability.

Recall that a subset $D$ of $X$ is a regular $G_\dlt$ if, for each $n\in\w$,
there are open sets $U_n$ containing $D$ such that $D=\bigcap_{n\in\w}\ol{U}_n$.
We say that a space is $\dlt$-normal if and only if every pair of
disjoint closed sets, one of which is a regular $G_\dlt$, can be separated
by disjoint open sets.
(In \cite{bu2} $\dlt$-normality refers to a distinct property.)

\begin{thm} For $T_3$ spaces the following are equivalent.
\begin{enumerate}
\item $X$ is countably paracompact.

\item For every decreasing sequence $(D_n)_{n \in \omega}$
of closed sets satisfying $\bigcap_{n\in\omega}D_n=\nowt$
there exists a sequence $(U_n)_{n \in \omega}$ of open sets
such that $D_n \subseteq U_n$ for $n \in \omega$ and
$\bigcap _{n \in \omega} \ol{U}_n =\nowt$ (Ishikawa~\cite{ish}).
\item For every countable increasing open cover $\{U_n \colon n \in \w\}$
there is an open refinement $\{V_n \colon n\in\w\}$ such that
$\ol{V}_n \subseteq U_n$ for all $n \in \w$.
\item For every decreasing sequence $(D_n)_{n \in \omega}$
of closed, nowhere dense sets satisfying $\bigcap_{n\in\omega}D_n=\nowt$
there exists a sequence $(U_n)_{n \in \omega}$ of open sets
such that $D_n\subseteq U_n$ for $n \in \omega$ and
$\bigcap _{n \in \omega} \ol{U}_n =\nowt$
(Hardy, Juhasz \cite{hj}).
\item For every increasing dense open cover there is a countable
closed refinement of sets whose interiors cover $X$.

\item For every  lower semicontinuous, real-valued function $g>0$
on $X$ there is an upper semicontinuous $h$ such that
$0<h<g$ (Dowker \cite{dowk}).
\item For every lower semicontinuous, real-valued function $g>0$
on $X$ there exist real-valued functions \(u\) and \(v\) where \(u\) is upper
semicontinuous and \(v\) is lower semicontinuous and such that
$0<u \les v \les g$
(Mack~\cite{mack2}).
\item For every locally bounded, real-valued function $g$ on $X$ there
exists a locally bounded, lower semicontinuous function $h$ such that
$|g| \les h$ (Mack~\cite{mack1}).

\item $X \times [0,1]$ is $\delta$-normal (Mack~\cite{mack2}).
\item If $C$ is a closed subset of $X \times [0,1]$ and $D$ is a closed
subset of $[0,1]$ such that $C \cap (X \times D) =\nowt$, then
there are disjoint open sets separating $C$ and $X\times D$ (Tamano~\cite{t}).
\end{enumerate}\label{cp equivs}
\end{thm}

One possible natural monotone version of countable paracompactness,
introduced independently in \cite{gks}, \cite{txl} and \cite{pan},
is MCP (that the definition stated below is equivalent to the
original is shown in \cite{gy}).

\begin{defn}
A space \(X\) is MCM if and only if there is an operator \(U\)
assigning to each \(n \in \omega\) and each closed set \(D\) an open
set \(U(n,D)\) such that
\begin{enumerate}
\item \(D \subseteq U(n, D)\),
\item if \(E \subseteq D\), then \(U(n, E) \subseteq U(n,D)\) and
\item if \((D_i)_{i \in \omega}\) is a decreasing sequence of closed
sets with empty intersection, then \(\bigcap_{n \in \omega} U(n,
D_n) =\n\).

\end{enumerate}

$X$ is MCP if, in addition, $U$ satisfies
\begin{enumerate}
\item[($3'$)] if \((D_i)_{i \in \omega}\) is a decreasing sequence of
closed sets with empty intersection,
then \(\bigcap_{n \in \omega} \ol{U(n, D_n)} =\n\).
\end{enumerate}
\end{defn}

MCM (monotonically countably metacompact) spaces are precisely the
$\bt$-spaces (see \cite{gks}) and MCP is closely related to
the wN property and stratifiability, which can be characterized (see \cite{gy})
by
conditions (1), (2) and
\begin{enumerate}
\item[(\(3''\))] if \((D_i)_{i \in \omega}\) is a decreasing sequence of closed
sets, then \newline \(\bigcap_{n \in \omega} \ol{U(n,D_n)} = \bigcap_{n \in
\omega} D_n\).
\end{enumerate}

It turns out that many of the natural monotone versions of the
properties listed in Theorem \ref{cp equivs} are equivalent to MCP.

\begin{thm} \label{mcp equivs}
For a $T_3$ space, the following are equivalent.
\begin{enumerate}
 \item \label{mcp} $X$ is MCP.

 \item \label{mcp old} There is an operator $V$ assigning to
each decreasing sequence $(D_j)_{j\in\omega}$ of closed sets with empty
intersection, a sequence of open sets
$(V(n,(D_j)))_{n\in\omega}$ such that:
$D_n\subseteq V(n,(D_j))$ for each $n\in\omega$;
if $D_n\subseteq E_n$ for each $n\in\w$, then $V(n,(D_j))\subseteq V(n,(E_j))$;
$\bigcap_{n\in\omega}\overline{V(n,(D_j))}=\varnothing$.

 \item \label{mcp pan}
 Suppose that $\mathbb{H}$ is any partially-ordered set and $F$ is any map
from $\omega\times\mathbb{H}$ to the closed subsets of $X$ such that
both $F(.,h)$ and $F(n,.)$ are order-reversing, and
$\bigcap_{n\in\omega}F(n,h)=\varnothing$ for all $h\in\mathbb{H}$.
There is a map $G$ from $\omega\times\mathbb{H}$ to the
open subsets of $X$ such that $F(n,h)\subseteq G(n,h)$, for each
$n\in\omega$ and $h\in\mathbb{H}$,
both $G(.,h)$ and $G(n,.)$ are order-reversing, and
$\bigcap_{n\in\omega}\overline{G(n,h)}=\varnothing$ for all
$h \in \mathbb{H}$.

 \item \label{mcp deM} There is an operator \(W\) assigning to each \(n \in
\omega\)
and each open set \(U\) an open set \(W(n,U)\) such that:
\(\ol{W(n,U)} \sq U\);
if \(U \sq U'\), then \(W(n,U) \sq W(n,U')\); and
if \((U_i)_{i \in \omega}\) is an increasing open cover of \(X\), then
\((W(n,U_n)^\circ)_{n \in \omega}\) is a refinement of
\((U_i)_{i \in \omega}\).

 \item \label{mcp dow} There is an operator $\varphi$ assigning to each
lower semicontinuous, real-valued function $0<g$ on $X$, an upper
semicontinuous $\varphi (g)$ such that: $0<\varphi (g)<g$; and
$\varphi(g)\leqslant \varphi(g')$ whenever $g\leqslant g'$.

\item \label{mcp mac0} There is an operator $\chi$ assigning to each
lower semicontinuous, real-valued function $0<g$ on $X$ and each
$i=0,1$, real-valued functions $\chi(g,0)$ and $\chi(g,1)$ such
that: $\chi(g,0)$ is lower semicontinuous, $\chi(g,1)$ is upper
semicontinuous; $0<\chi(g,0) \les \chi(g,1) \les g$ (equivalently $0<\chi(g,0) <
\chi(g,1) < g$); and $\chi(g,i)\les
\chi(g',i)$, $i=0,1$, whenever $g\les g'$.

 \item \label{mcp mac} There is an operator $\psi$ assigning to each
locally bounded, real-valued function $g$ on $X$, a locally bounded,
lower semicontinuous, real-valued $\psi(g)$ such that: $|g|\les
\psi(g)$ (equivalently $|g|<\psi(g)$); and $\psi(g)\les \psi(g')$
whenever $|g|\les|g'|$.

 \item \label{mcp mac3} There is an operator $\pi$ assigning to each
locally bounded, real-valued function $g$ on $X$ and each $i=0,1$,
locally bounded, real-valued functions $\pi(g,0)$ and $\pi(g,1)$
such that; $\pi(g,1)$ is lower semicontinuous, $\pi(g,0)$ is upper
semicontinuous; $\pi(g,0)\les g\les \pi(g,1)$ (equivalently
$\pi(g,0)< g<\pi(g,1)$); and $\pi(g,i)\les \pi(g',i)$, $i=0,1$,
whenever $g\les g'$.

 \item \label{mcp prod} $X\times[0,1]$ is MCP.

 \item \label{mcp tam} There is an operator \(H\) assigning to each
pair \((C,D)\), where \(C\)  is closed in \(X \times [0,1]\) and
\(D\) is closed in \([0,1]\) such that \(C\cap(X\times D)=\n\),
an open set \(H(C,D)\) such that:
 \begin{enumerate}
 \item \(C \sq H(C,D) \sq \ol{H(C,D)} \sq X \times ([0,1]\sm D)\);
 \item if \(C \sq C'\) and \(D' \sq D\), then \(H(C,D) \sq H(C',D')\).
 \end{enumerate}
\end{enumerate}
\end{thm}

\begin{pf} The equivalences of (\ref{mcp}), (\ref{mcp old}),
(\ref{mcp pan}), (\ref{mcp dow}),  (\ref{mcp prod}), and (\ref{mcp
mac} for $\les$) are proved in \cite{gks} and \cite{gy}. That (\ref{mcp
deM}) is equivalent to (\ref{mcp}) follows easily by de Morgan's
Laws. Clearly (\ref{mcp mac0}, \ref{mcp mac}, \ref{mcp mac3} for $<$) imply
(\ref{mcp mac0}, \ref{mcp mac}, \ref{mcp mac3} for $\les$).

(\ref{mcp mac0}) \(\rightarrow\) (\ref{mcp mac} for $<$):
Suppose \(g \colon X \to \mathbb R\) is locally bounded.
Define \(h = 1/(|g| +1)^*\)
(where \(f^*\) denotes the upper limit of \(f\)).  Then \(h\)
is lower semicontinuous and strictly positive.  By (6) there
exists a real-valued upper semicontinuous function \(\chi(h,1)\)
such that \(0 < \chi(h,1) \les h\).  Let \(\psi(g) = 1/\chi(h,1)\).
Then \(\psi(g)\) is lower semicontinuous and \(|g| < (|g| +1)^* \les \psi(g)\).
It is routine to show that \(\psi(g) \les \psi(g')\) whenever
\(|g| \les |g'|\).

(\ref{mcp mac}) \(\rightarrow\) (\ref{mcp mac3} for $<$):
If $g$ is a
locally bounded function, then $g^+=\max\{g,0\}$,
$g^{-}=\min\{g,0\}$ and $-g^{-}$ are locally bounded and $g^{-} \les
g\les g^{+}$. By (\ref{mcp mac}), there are locally bounded, lower
semicontinuous functions $\psi(g^{+})$ and $\psi(-g^{-})$ such that
$|g^{+}|=g^{+} < \psi(g^{+})$ and $|g^{-}|=-g^{-} < \psi(-g^{-})$. But then
$-\psi(-g^{-}) < g^{-} \les g\les g^{+}
< \psi(g^{+})$, $\psi(g^{+})$ is locally bounded and lower
semicontinuous, and $-\psi(-g^{-})$ is locally bounded and upper
semicontinuous. Moreover, if $g\les g'$, then $g^{+} \les (g')^{+}$,
so that $\psi(g^{+})\les \psi((g')^{+})$, and $0\leqslant -(g')^{-}
\les -g^{-}$, so that $-\psi(-g^{-}) \les -\psi(-(g')^{-})$.

(\ref{mcp mac3}) \(\rightarrow\) (\ref{mcp mac}):
Suppose \(g \colon X \to \mathbb R\) is locally bounded.
Define \(h = (|g| +1)^*\).  Then \(h\) is upper semicontinuous
and strictly positive and is therefore locally bounded.  By (8),
there exists a locally bounded, lower semicontinuous function \(\pi(h,1)\)
such that \(h \les \pi(h,1)\).  Let \(\psi(g) = \pi(h,1)\), then
\(|g| < h \les \psi(g)\).  If \(|g| \les |g'|\), then \(h \les h'\)
and so \(\psi(g) \les \psi(g')\).

(\ref{mcp tam}) $\rightarrow$ (\ref{mcp mac0} for $<$): Suppose that
$0<g$, where $g$ is a lower semicontinuous function. We follow the
argument used in the proof of Theorem 1.6 of \cite{gs}.
$B_g=\{ (x,r) \colon r\geqslant g(x)\}$ and $X\times \{0\}$ are disjoint
closed subsets of $X\times [0,1]$, so by (\ref{mcp tam}), there exists
an open set $V_g$ such that $X\times\{0\}\subseteq V_g\subseteq
\ol{V}_g\subseteq \big(X\times[0,1]\big)\smallsetminus B_g$ with the
property that $V_g\subseteq V_{g'}$ whenever $g\les  g'$. Let
\begin{align*}
u_g(x)&=\sup\{ r\colon (x,s)\in\overline{V}_g\textrm{ for
all }s<r\}\\
l_g(x)&=\sup\{ r\colon (x,s)\in V_g\textrm{ for
all }s<r\}.
\end{align*}
As in \cite{gs}, $u_g<g$, $u_g\les u_{g'}$, whenever $g\les g'$,
and $u_g$ is upper semicontinuous.
Similarly $0<l_g\les u_g$ and $l_g\les l_{g'}$, whenever
$g\les g'$. Monotonicity is also clear.

It remains to prove that $l_g$ is lower semicontinuous,
to which end we show that
$l^{-1}_g(r,\infty)$  is open for any $r\in \R$.
Suppose, then that $x\in l^{-1}_g(r,\infty)$ so that $r<l_g(x)$.
Choose some $r<t<l_g(x)$, then for each
$s\in[0,t]$, $(x,s)\in V_g$. For each such $s$, choose an
open neighbourhood $W_s$ of $x$ and an $\eps_s>0$ such that
$W_s\times(s-\eps_s,s+\eps_s)\subseteq V_g$. By compactness,
there are finitely many $s_0,\dots,s_n$ such that
$\bigcup_{i\les n}\big(W_{s_i}\times(s_i-\eps_{s_i},s_i+\eps_{s_i})\big)$
covers $\{x\}\times[0,t]$. Hence
$(x,r)\in\big(\bigcap_{i\les n} W_{s_i}\big)\times(0,t)$, so for all
$y\in\bigcap_{i\les n} W_{s_i}$, $l_g(y)>r$ as required.


(\ref{mcp}) \(\rightarrow\) (\ref{mcp tam}): Finally, let $U$ be an MCP
operator on $X$ such that $U(n+1,E)\subseteq U(n,E)$ for each closed $E$
and $U(n,\nowt)=\nowt$.
Suppose that \(D\) is any closed subset of \([0,1]\) and
$C$ is any closed subset of $X\times[0,1]$ disjoint from $X\times D$.

Let \(V(n,D) = \{s \in [0,1] \colon (\exists d\in D)|s-d|<1/n\}\) and let
$C(n,D)$ be the projection onto $X$ of $C \cap \ol{X \times V(n,D)}$.
Then $(V(n,D))_{n \in \omega}$ is a decreasing sequence of open sets
containing $D$ such that \(D= \bigcap _{n \in \omega} \ol{V(n,D)}\)
and \((C(n,D))_{n\in\w}\) is a decreasing sequence of closed subsets of
$X$ with empty intersection. Furthermore,
\(\bigcap _{n \in \omega} \ol{U(n, C(n,D))} =\n\)
and \(U(m+1,C(n,D)) \sq U(m,C(n,D))\) for each \(m \in
\omega\).

Now, for $r\in[0,1]$, let \(\varepsilon _{r,D} = \inf\{|r-d| \colon d\in D\}/2\)
and,
if $r\notin D$, let \(n_{r,D}\) be the least natural number such that
\(1/n_{r,D}<\varepsilon_{r,D}\). Let \(C_r\) be the closed set
$\{x\in X \colon (x,r)\in C \cap (X\times \{r\})\}$. Then \(C_r \sq U(n_{r,D},
C_r)\) and
\(U(n_{r,D}, C_r) \sq U(n_{r,D},C_r')\), whenever \(C_r \sq C_r'\).

Define \(H(C,D) = \bigcup _{r \in [0,1]} U(n_{r,D}, C_r) \times
(r-\varepsilon_{r,D}, r+\varepsilon _{r,D})\).
Clearly $H(C,D)$ is open and contains \(C\).
Suppose \(C \sq C'\) and \(D' \sq D\).  For any
\(r\notin D\), \(n_{r,D'} \les n_{r,D}\).  Hence, by monotonicity,
\(U(n_{r,D}, C_r) \sq U(n_{r,D'},
C_r')\), from which it follows that \(H(C,D) \sq H(C',D')\).

It remains to show that \(\ol{H(C,D)}\) and \(X\times D\) are disjoint.
To this end, let \((x,d) \in X \times D\).
From above, there is some \(n_x\in \omega\) such that
$x$ is not in $\ol{U(n_x, C(n_x,D))}$. Let $W$ be an
open neighbourhood of \(x\) disjoint from $U(n_x, C(n_x,D))$.
For any \(r \in V(n_x,D)\), \(n_{r,D} > 2n_x\) and
\(C_r \sq C(n_x,D)\).
Hence \(U(n_{r,D}, C_r) \sq U(n_{r,D},C(n_x,D))\sq U(n_x,C(n_x,D)))\)
and $W$ is disjoint from $U(n_{r,D}, C_r)$ for any \(r \in V(n_x,D)\).

Now, if $r\notin V(n_x,D)$, then $J_d=\bigg(d-\frac1{2n_x},
d+\frac1{2n_x}\bigg)$ and
$(r-\varepsilon _{r,D}, r+\varepsilon _{r,D})$ are disjoint, so that
\[(W \times J_d) \cap \bigg(U(n_{r,D}, C_r) \times
(r-\varepsilon _{r,D}, r+\varepsilon _{r,D}) \bigg) = \n\]
for any \(r \in [0,1]\).  Thus
\(W \times J_d\) is an open
neighbourhood of \((x,d)\), disjoint from \(U(C,D)\).
\end{pf}

We note that essentially the same proof shows that conditions (9)
and (10) hold for any compact metric space in place of $[0,1]$.

The natural monotone version of $\dlt$-normality seems to be the following.

\begin{defn}
A space $X$ is \textit{monotonically} \(\delta\)\textit{-normal}
($\mdn$) iff there is an operator $H$ assigning to each pair
of disjoint closed sets $C$ and $D$, at least one of which is a regular
\(G_{\delta}\)-set, an open set $H(C,D)$ such that
\begin{enumerate}
\item[(1)] $C\subseteq H(C,D)\subseteq \ol{H(C,D)}\subseteq X\setminus D$
\item[(2)] if $C\subseteq C'$ and $D'\subseteq D$, then
$H(C,D)\subseteq H(C',D')$.
\end{enumerate}
\end{defn}

Replacing \(H(C,D)\) with \(H(C,D) \sm \ol{H(D,C)}\), if necessary,
one may assume that \(H(C,D) \cap H(D,C)=\n\).

We discuss $\mdn$ in more detail in \cite{gh}, for the present we restrict
our attention to the relationship between MCP and $\mdn$.

\begin{thm}\label{thm5}
If \(X \times [0,1]\) is $\mdn$, then \(X\) is both $\mdn$  and MCP.
Moreover, if $X\times [0,1]$ is $\mdn$, then $X\times [0,1]$ is MCP.
\end{thm}

\begin{pf} If $C$ and $D$ are disjoint closed subsets of
$X$, one of which is a regular $G_\dlt$, then $C\times \{0\}$
and $D\times \{0\}$ are disjoint closed subsets of $X\times [0,1]$,
one of which is a regular $G_\dlt$. Furthermore, if $E$ is
a closed subset of $[0,1]$, then $X\times E$ is a regular $G_\dlt$
subset of $X\times [0,1]$.
It follows that, if $X\times [0,1]$ is $\mdn$, then $X$ is $\mdn$ and
satisfies condition (\ref{mcp tam}) of Theorem \ref{mcp equivs}.
The last statement
of the theorem follows by condition (\ref{mcp prod}) of Theorem
\ref{mcp equivs}.
\end{pf}

\begin{exmp}
There exists an MCP space that is not $\mdn$,  hence its product
with \([0,1]\) is not $\mdn$.
\end{exmp}

\begin{pf} Let $\A$ be a maximal almost disjoint family of subsets of
$\w$, then Mr\'owka's $\Psi$ \cite{mro} is the locally compact, locally
countable, zero-dimensional, non-metrizable Moore space,
$\Psi=\w\cup\A$, in which each $n\in\w$ is isolated and basic
(clopen) neighbourhoods of $a\in \A$ take the form
$\{a\}\cup\big(a\cap(k,\w)\big)$
for some $k\in \w$. Since $\Psi$ is a non-metrizable Moore space,
it is not monotonically normal.

Let $X$ be the one-point compactification of $\Psi$, so that $X$
is MCP (see \cite{gks}). We claim that $X$ is not $\mdn$.

Observe that if $C$ and $D$ are disjoint closed subsets
of $X$, at least one contains at most a finite subset of $\A$.
It is easy to show that any closed subset of $X$ containing at most finitely
many points of $\A$ is a regular \(G_{\delta}\) subset of $X$.
It follows that, if $X$ were $\mdn$, it would be monotonically normal,
and hence $\Psi$ would be monotonically normal.

By Theorem \ref{thm5}, the product $X\times [0,1]$ is not $\mdn$.
\end{pf}

\begin{exmp}\label{sorgen}
There exists an $\mdn$ space $X$ such that $X$ is not MCP
and $X\times [0,1]$ is not $\mdn$.
\end{exmp}

\begin{pf} It is shown in \cite{gks} that both the Sorgenfrey and
Michael lines are monotonically normal, so $\mdn$, but not MCP.
\end{pf}

\begin{exmp} There is a space $X$ that is both $\mdn$ and MCP
such that $X\times[0,1]$ is not $\mdn$.
\end{exmp}

\begin{pf}
The Alexandroff duplicate of the unit interval $[0,1]$ is compact,
therefore MCP, and monotonically normal but not stratifiable (as it is not
perfect).
However, it is first countable and regular, so by a result in
\cite{gh}, if $X\times [0,1]$ were $\mdn$, it would be monotonically normal and
therefore
$X$ would be stratifiable.
\end{pf}

The appropriate monotonization of Theorem \ref{cp equivs} (4) is clearly
the following.

\begin{defn} A space $X$ is nMCM (nowhere dense MCM)
if there is an operator $U$ assigning an open set $U(n,D)$
to each natural number $n\in\w$ and each closed, nowhere dense set $D$
such that:
\begin{enumerate}
\item $D\sbs U(n,D)$,
\item if $E\subseteq D$ then $U(n,E)\subseteq U(n,D)$ and
\item if $(D_n)$ is a decreasing sequence of closed,
nowhere dense sets with
empty intersection, then $\bigcap U(n,D_n)=\nowt$.
\end{enumerate}
If, in addition, $U$ satisfies
\begin{enumerate}
\item[($3^\prime$)] if $(D_n)$ is a decreasing sequence of closed,
nowhere dense sets with
empty intersection, then $\bigcap\ol{U(n,D_n)}=\nowt$,
\end{enumerate}
then $X$ is said to be nMCP.
\end{defn}

It turns out that nMCM is equivalent to MCM, but that the situation for nMCP is
not clear.

Let $X$
be a space and, for each $x\in X$ and $n\in\omega$ let $g(n,x)$ be an
open set containing $x$. We say that $X$ is a $\beta$-space, or
wN-space, or q-space
(see \cite{hodel1,hodel2,michael}), if it satisfies
\begin{enumerate}
\item[($\beta$)] if $x\in g(n,y_n)$ for all $n$, then the sequence
$(y_n)$ has a cluster point;
\item[(wN)] if $g(n,x)\cap g(n,y_n)\neq\varnothing$ for all $n$, then
the sequence $(y_n)$ has a cluster point;
\item[(q)] if $y_n\in g(n,x)$ for all $n$, then the sequence $(y_n)$
has a cluster point.
\end{enumerate}

\begin{thm} \begin{enumerate}
\item $X$ is nMCM if and only if it is MCM
if and only if it is a $\beta$-space.
\item Suppose that $X$ is either a $q$-space or a locally
countably compact space. $X$ is nMCP if and only if it is
MCP if and only it is a wN space. In particular every
first countable or locally compact
nMCP space is MCP.
\end{enumerate}
\end{thm}

\begin{pf} For (1) recall that $\beta$-spaces are
precisely MCM spaces
\cite{gks}. Clearly every MCM space is nMCM.

We modify the proof used in \cite{gks} to show that every
nMCM space is a $\bt$-space. To this end, for each $x$ in
$X$ and $n$ in $\omega$ let
$D^n_j(x)=\{x\}$ if $j\leqslant n$ and
$D^n_j(x)=\varnothing$ otherwise. For each $x\in X$ and fixed $n$,
$(D^n_j(x))_{j\in\omega}$ is a decreasing sequence of closed sets with
empty intersection.  Now let
$$
g(n,x)=
\begin{cases}
\{x\} & \mbox{if $x$ is isolated}\\
U(n,D^n_j(x)) & \mbox{if $x$ is not isolated.}\\
\end{cases}
$$ Suppose that
$(y_n)$ is a sequence of distinct points without a cluster point.
There are two cases: either infinitely many $y_n$ are isolated,
or at most finitely many are isolated. In the first case,
for some $n$ and $m$, $g(n,y_n)\cap g(m,y_m)=\nowt$ and so
$\bigcap_{n\in\omega}g(n,y_n)=\varnothing$. In the second case,
we may, without loss of
generality, assume that no $y_n$ is isolated and
define $E_j=\{y_n \colon n\geqslant j\}$, for each $j\geqslant k$.
Clearly $(E_j)$ is a decreasing sequence of closed, nowhere dense
sets with empty intersection and for all $n,j\in\omega$,
$D^n_j(y_n)\subseteq E_j$. By monotonicity $g(n,y_n)\subseteq
U(n,E_n)$ for each $n$ and therefore
$\bigcap_{n\in\omega}g(n,y_n)=\varnothing$. Thus $X$ is a
$\beta$-space.

The proof of (2) is similar. Define $g$ as in (1).
If $X$ is a q-space, then, for each $x\in X$ and
$n\in\omega$, choose $h(n,x)$ satisfying condition (q) above; if $X$
is locally countably compact then let $h(n,x)=C_x$ for each $x$ and
$n$ where $C_x$ is a neighbourhood of $x$ with countably compact
closure.  Let $G(n,x)=g(n,x)\cap h(n,x)$. We claim that the $G(n,x)$
satisfy the conditions for a wN-space. Clearly $x\in G(n,x)$ for each
$x$ and $n$.  So assume that $x_n\in G(n,x)\cap G(n,y_n)$ for each
$n$.

If $y_n$ is isolated, then $G(n,y_n)=\{y_n\}$ and
$y_n\in G(n, x)$. Hence, if infinitely many $y_n$ are isolated, then
either by q or by local countable compactness, they cluster.
So suppose that at most finitely many $y_n$ are isolated.
Again we may assume in fact that none are isolated.

Either by the condition (q), or by the countable compactness of
$\ol{C}_x$, the sequence $(x_n)$ has a cluster point $z$. Assume
for a contradiction that the sequence $(y_n)$ does not have a cluster
point. Define $E_j$ as in (1).
Again
$(E_n)$ is a decreasing sequence of closed, nowhere dense
sets with empty
intersection and, by nMCP, $x_n\in G(n,y_n)\subseteq U(n, E_n)$ for
each $n$. Now without loss of generality we may assume that
$U(m,(E_j))\subseteq U(n,(E_j))$ for all $m\geqslant n$ and so, for
each $n$, $x_m\in U(n,(E_j))$ for all $m\geqslant n$. Since the $x_n$
cluster at $z$, this implies that $z\in \overline{U(n,(E_j))}$ for each
$n$ which is a contradiction.
\end{pf}

Similar arguments show that many results about MCP spaces hold for
nMCP spaces,
for example: every nMCP, Moore space is metrizable;
first countable nMCP spaces are collectionwise Hausdorff;
if a nMCP space is not collectionwise Hausdorff, then there is a
measurable cardinal. However, we have been unable to determine whether nMCP and
MCP coincide.

\medskip

Finally let us note that Dowker \cite{dowk} proved the following theorem.

\begin{thm} \label{dowker3}
The following properties of a Hausdorff space \(X\) are equivalent:
\begin{enumerate}
\item \(X\) is countably paracompact and normal.
\item \(X\) is countably metacompact and normal.
\item If \(g \colon X \to \mathbb R\) is lower semicontinuous,
\(h \colon X \to \mathbb R\) is upper semicontinuous,
and \(h < g\), then there exists a continuous
function \(f \colon X \to \mathbb R\) such that \(h < f < g\).
\item \(X \times [0,1]\) is normal.
\end{enumerate}
\end{thm}

The monotone versions of these conditions are not equivalent \cite{gks,gs}. A
space is
monotonically normal and MCM iff it is monotonically normal and MCP. On the
other hand,
a space is stratifiable iff whenever \(g \colon X \to \mathbb R\) is lower
semicontinuous
and \(h \colon X \to \mathbb R\) is upper semicontinuous
and \(h < g\), then there exists a continuous
function \(f(g,h)\) such that \(h < f(g,h) < g\) such that
$f(g,h)\les f(g',h')$ whenever $g\les g'$ and $h\les h'$.


\begin{thebibliography}{99}
\bibitem{as}R. A. Al\`{o} and H. L. Shapiro, \textit{Normal topological
spaces}, Cambridge University Press, 1974.
\bibitem{bu} D. Burke, \textit{PMEA and first countable, countably paracompact
spaces}, Proc. Amer. Math. Soc. 92 (1984), no. 3, 455--460.
\bibitem{bu2} D. Burke, \textit{Covering properties}, Handbook of
Set-Theoretic Topology, K. Kunen and J. E. Vaughan, North-Holland, Amsterdam
(1984), 347--422.
\bibitem{dowk}C. H. Dowker, \textit{On Countably Paracompact Spaces}, Canad. J.
Math. 3 (1951), 219--224.
\bibitem{eng} R. Engelking, \textit{General topology}, Heldermann Verlag,
Berlin,
1989.
\bibitem{gh} C. Good and L. Haynes, \textit{Monotone $\delta$-normality},
in preparation.
\bibitem{gk} C. Good and R. Knight, \textit{Monotonically countably paracompact,
collectionwise Hausdorff spaces and measurable cardinals},  Proc. Amer. Math.
Soc. 134 (2006), no. 2, 591--597.
\bibitem{gks} C. Good, R. Knight and I. Stares, \textit{Monotone countable
paracompactness}, Topology Appl. 101 (2000), no. 3, 281--298.
\bibitem{gs} C. Good and I. Stares, \textit{Monotone insertion of continuous
functions}, Topology Appl. 108 (2000), no. 1, 91--104.
\bibitem{gy} G. Ying and C. Good, \textit{A note on monotone countable
paracompactness}, Comment. Math. Univ. Carolin. 42 (2001), no. 4, 771--778.
\bibitem{hj} K. Hardy and I. Juh\'asz, \textit{Normality and the weak cb
property}, Pacific J. Math. 64 (1976), no. 1, 167--172.
\bibitem{hodel1} R. E. Hodel, \textit{Moore spaces and $w\Delta$-spaces},
Pacific J. Math. 38 (1971), 641--652.
\bibitem{hodel2} R. E. Hodel, \textit{Spaces defined by sequences of open covers
which guarantee that certain sequences have cluster points}, Duke Math. J. 39
(1972), 253--263.
\bibitem{ish} F. Ishikawa, \textit{On countably paracompact spaces}, Proc.
Japan Acad. 31 (1955), 686--687.
\bibitem{kat}
M. Kat\v{e}tov, \textit{Measures in fully normal spaces}, Fund.
Math. 38 (1951), 73--84.
\bibitem{hand} K. Kunen and J. E. Vaughan, eds., Handbook of
set-theoretic topology, (North-Holland, Amsterdam, 1984).
\bibitem{mack1} J. E. Mack, \textit{On a class of countably paracompact spaces},
Proc. Amer. Math. Soc. 16 (1965), 467--472.
\bibitem{mack2} J. E. Mack, \textit{Countable paracompactness and weak normality
properties}, Trans. Amer. Math. Soc. 148 (1970), 265--272.
\bibitem{michael} E. Michael, \textit{A note on closed maps and compact sets},
Israel J. Math. 2 (1964), 173--176.
\bibitem{mro} S. Mr\'{o}wka, \textit{On completely regular spaces}, Fund. Math.
41 (1954), 105--106.
\bibitem{pan} C. Pan, \textit{Monotonically CP spaces}, Q \& A in General
Topology 15 (1997), no. 1, 25--32.
\bibitem{t} H. Tamano, \textit{On compactifications}, J. Math. Kyoto Univ. 1
(1962), 161--193.
\bibitem{txl} H. Teng, S. Xia and S. Lin, \emph{Closed images of some
generalized countably
compact spaces}, Chinese Ann. Math. Ser A, 10 (1989),
554--558.
\bibitem{wat} W. S. Watson, \textit{Separation in countably paracompact spaces},
Trans. Amer. Math. Soc. 290 (1985), 831--842.
\end{thebibliography}
\end{document}